\documentclass[12pt]{amsart}
\usepackage{amssymb}
\newcommand{\pp}{{{\Bbb P}}}

\newcommand {\OO} {\mathcal O}
\newcommand {\II} {\mathcal I}
\newcommand {\ZZ}{{\Bbb{Z}}}

\newcommand\Pic{\operatorname{Pic}}
\newtheorem{lemma}{Lemma}[section]
\newtheorem{prop}[lemma]{Proposition}
\newtheorem{thm}[lemma]{Theorem}
\newtheorem{corol}[lemma]{Corollary}
\newtheorem{rem}[lemma]{Remark}
\newtheorem{defn}[lemma]{Definition}
\newtheorem{ex}[lemma]{Example}

\begin{document}
\title{Halphen conditions and postulation of nodes}

\author[L.\ Chiantini, N.\ Chiarli, S.\ Greco]
{L.\ Chiantini$^1$, N.\ Chiarli$^2$, S.\ Greco$^3$}
\address{$^1$Dipartimento di Matematica, Universit\`a di Siena,  I-53100
Siena, Italy}
\email{CHIANTINI@UNISI.IT}
\address{$^2$ Dipartimento di Matematica, Politecnico di Torino,  I-10129
Torino, Italy}
\email{NADIA.CHIARLI@POLITO.IT}
\address{$^3$ Dipartimento di Matematica, Politecnico di Torino,  I-10129
Torino, Italy}
\email{SILVIO.GRECO@POLITO.IT}

\subjclass{14H50}

\thanks{Supported by MIUR and GNSAGA-INDAM, Italy}


\begin{abstract} We give sharp lower bounds for the postulation of the
nodes of a general plane projection
of a smooth connected curve $C \subseteq \pp^r$ and we study the
relationships with the geometry of the
embedding. Strict connections with Castelnuovo's theory and Halphen's
theory are shown.
\end{abstract}

\bigskip

\date{\today}

\maketitle
\tableofcontents


\section*{Introduction}

This paper deals with smooth projective non-degenerate curves $C\subset
\pp^r$, of degree $d$ and genus $g$, and their general projections $C\to
C_0\subset \pp^2$. If $N$ is the set of nodes of $C_0$, it is classically known that the
postulation of $N$ determines many intrinsic and extrinsic geometric
properties of $C$, such as, for example, the existence of special linear series, the classical relation among $d$, $g$ and the degree $\delta$ of $N$ and the equality
$e = d-3-\alpha$, where $e$ is the speciality of $C$ and  $\alpha$ is the
least degree of a curve containing $N$ (see Section 1).

Much subtler results show how the existence of special linear series
$g^n_k$ depends on the existence of surfaces of low degree $s$ containing
$C$, via complicate but sometimes sharp numerical relations between $r$,
$n$, $k$, $s$, $d$ and $g$, that one can find in the literature (see \cite {A}, \cite {AS}, \cite {CCD1}, \cite{CL}, \cite{CK}, or \cite{CC}, in which our principal procedure is settled). \par

The main goal of this paper is to prove some sharp lower bounds for
$\partial h_N$, the first difference of the Hilbert function of $N$, and
their relationships with the geometry of the embedding $C\hookrightarrow
\pp^r$. Since $\alpha$ is also the maximum of $\partial h_N$, this approach brings
into play,  in a natural way, the speciality $e$. Indeed our bounds will be
given in terms of $d$, $r$ and $\alpha$ (rather than $e$) and, when $r =
3$, we consider also an integer $s$ such
that $C$ is contained in no surface of degree $s-1$. These bounds are
obviously related to upper bounds for $g$ in terms of $d$, $r$, $\alpha$
(hence $e$) and $s$. So our theory is strictly related to Castelnuovo's
theory and to Halphen's theory of maximal genus, on which we
hope to bring some new point of view.

The paper is organized as follows. After the preliminary Section 1, in the
second section, we deal with non-degenerate curves in $\pp^r$, for arbitrary $r\ge3$ and we give a sharp lower bound for $\partial h_N$, depending only on degree and
speciality, using the main result in \cite {CGN}. The result is a Castelnuovo-like theory on the relationship between postulation of nodes and genus of curves.

In the rest of the paper we work in $\pp^3$ and we bring $s$ into play. In
Section 3 we give some general results which show how  $s$ affects
$\partial h_N$, the most used in the sequel being Theorem \ref {p3bis}. For
this we use two main tools: first we translate the Hilbert function of $N$
into the dimensions of some linear series on $C$, using the classical
approach and the results contained in \cite{CGN}; next we bound the dimensions of
the linear series, using the procedure introduced in \cite{CC} and the
results of \cite{CCD2}.

In Section 4 we prove the key Theorem \ref{arit}, which produces our main
lower bound for $\partial h_N$, depending on $d$, $s$ and $\alpha$.  As
immediate but useful Corollaries we get an upper bound for $g$ (in terms of $d$, $s$ and
$\alpha$) and a lower bound for $\alpha$ (in terms of $d$ and $s$).

Then we study the sharpness of our bounds, with the extra assumption $d >
s(s-1)$. In this case our lower bound for $\alpha$ can be made explicit:
thus we get a purely elementary arithmetical proof of the celebrated
Halphen-Gruson-Peskine bound for $e$ (Proposition \ref{HGP}) and we can
show that the curves with $\alpha$ minimal and minimal node function are
exactly the Halphen curves. Indeed our approach allows an elementary proof
of some known characterizations of such curves (Theorem \ref {Halphen}).

The lower bounds for $\partial h_N$ is obtained so far when the value of
$\alpha$ is minimal and the curves attaining the bound turn out to be aCM.

Our next step, developed in Section 5, is to consider minimal nodal
functions with $\alpha$ non-minimal. Here the geometry becomes much more
difficult, and we are able to give some results only for $s \le 4$. The
case $s = 2$ is easy and we can give a complete picture. For $s=3$ and $s =
4$ we show that our bounds are sharp for every $d > s(s-1)$, by producing
explicitly curves on the smooth cubic surface or in the so-called Mori
surfaces (see \cite {M}). Moreover, using Theorem \ref{p3bis}, we show that, under some mild numerical restriction, every curve attaining the bound must lie on a surface of degree $s$. Surprisigly the curves we produce, except for some sporadic cases, are not
aCM: this makes us believe that the study of curves with minimal nodal postulation but non-minimal $\alpha$, and of their cohomology, is a challanging problem, even for low
values of $s$.

\section{Preliminaries and known results}\label{prel}

In this section we fix some notation and we collect some known results.

We work over an algebraically closed  field of characteristic zero.\par

Let $C\subset \pp^r$, $r \ge 3$, be a smooth connected non-degenerate curve
of degree $d$ and genus $g$.\par

We put $n:= h^0(\OO_C(1))-1$, whence $C$ is a projection of a linearly
normal curve of degree $d$ in $\pp^n$ (i.e. of a curve in $\pp^n$ which is
not a projection of a curve of the same
degree spanning a higher dimensional projective space).
In particular $C$ is linearly normal if and only if $n = r$.\par

We denote by $e$ the index of speciality of $C$, namely
$$
e:= \max\{j \mid h^1(\OO_{C}(j))>0\} = \max\{j \mid h^2(\II_{C}(j))>0\}.
$$
Recall that $C$ is said to be {\it special} if and only if $e > 0$.

Let $C\to C_0\subset \pp^2$ be a general plane projection, so that $C_0$
has only nodes as singularities. \par

Let $N \subset \pp^2$ be the reduced subscheme of the nodes of $C_0$ and
denote by $\Delta$ the pull-back of $N$ to $C$, viewed as a divisor on
$C$.\par

If $\delta$ is the length of $N$, we have the classical relation:
$$
\delta = \frac{(d-1)(d-2)}2 - g.
$$

Let $h_N$ be the Hilbert function of $N$ and let $\partial h_N$ be its
first difference.\par
Finally let
$$
\alpha:= \min \{j \; \vert\;  h^0(\II_N(j))\neq 0\}.
$$
\medskip

The classical theory of adjoint curves relates the linear systems of curves
passing through $N$ to the geometry of $C$. The starting points of this
theory can be summarized in the following remark (for details see e.g.
\cite{CGN}):

\begin{rem} \label{adjoints}\rm Let $D$ be a hyperplane divisor of $C$.
Then we have:
\begin{itemize}
\item[(i)] the canonical map $H^0( \II_N(j)) \to H^0(\OO_{C}(jD-\Delta))$
is surjective for every $j \in \ZZ$ and also injective for $j \le d-1$;
\item[(ii)] the canonical sheaf $\omega_{C}$ of $C$ is
$\OO_{C}((d-3)D-\Delta)$;
\item[(iii)] from (i) and (ii) it follows, by Riemann-Roch and standard
calculations:
$$
\partial h_N(j)=j+1-d-h^0(\OO_{C}(d-j-3))+h^0(\OO_{C}(d-j-2)) \quad {\rm
for}\; j \le d-1;
$$
\item[(iv)] by (ii) the the map of (i) can be viewed as $H^0( \II_N(j)) \to
H^0(\omega_{C}(d-3-j))$, hence from (iii) it follows:
$$
\partial h_N(j)=j+1-h^0(\omega_{C}(j-d+3))+h^0(\omega_{C}(j-d+2)) \quad
{\rm for}\; j \le d-1;
$$
\item[(v)] from the above we get $\alpha = d-3-e$. It follows $\alpha \le
d-2$, with equality if and only if $g=0$. Moreover $\alpha \le d-4$ if and
only if $C$ is special.
\end{itemize}

\end{rem}
\bigskip

By general properties of the Hilbert function of a zero-dimensional
subscheme of $\pp^2$ and
by Remark \ref{adjoints} we have the following result:
\medskip

\begin {prop}\label{basic} The difference function $\partial h_N$
satisfies:\par
\begin{itemize}
\item[(a)] $\partial h_N(j) = 0$ for $j < 0$ and for $j > d-3$;
\item[(b)] $\partial h_N(j) = j+1$ for $0 \le j \le\alpha-1$;
\item[(c)] $\partial h_N$ is non-increasing for $j\geq \alpha - 1$;
\item[(d)] $\partial h_N(d-3)=n-2 \ge r-2$.
\end{itemize}
\end{prop}
\medskip

Proposition \ref{basic}  above gives an accurate description of the
function $\partial h_N$, except for the interval $[\alpha,d-4]$.

The following result (\cite{CGN}, Theorem 5.1), which holds in arbitrary
characteristic,  provides a lower bound for the descent of $\partial
h_N(j)$ at any step between $\alpha$ and $d-3$. More refined bounds for $r
= 3$ will be given in the next sections.
\medskip

\begin{thm}\label{pr} {\rm (arbitrary characteristic)} For $\alpha\leq
j\leq d-3$ one has:
$$
\partial h_N(j+1)\leq \partial h_N(j)-n+2 \le \partial h_N(j)-r+2 .
$$
In particular $\partial h_N$ is strictly decreasing in the interval
$[\alpha, d-2]$.
\end{thm}

It is easy to derive some (well known) speciality bounds from the previous
result.

\begin{corol} \label{bound for alpha} The following inequalities hold:
\begin{itemize}
\item[(a)] $\alpha\geq  (d-2-\alpha)(n-2) \geq  (d-2-\alpha)(r-2);$
\item[(b)] $\alpha\geq  \frac{(d-2)(n-2)}{n-1} \geq  \frac{(d-2)(r-2)}{r-1};$
\item[(c)] $e \leq\frac {d-n-1}{n-1} \leq\frac {d-r-1}{r-1}.$
\end{itemize}
\end{corol}

\begin {proof} $\alpha$ is the maximum for $\partial h_N$. On the other
hand by Theorem \ref{pr} $\partial h_N$ decreases of at least $n-2$ at any
step from $j=\alpha$ to $j=d-3$; furthermore
$\partial h_N(d-3)=n-2$ by Proposition \ref {basic}(d) and (a) follows.

Clearly (b) is equivalent to (a). Finally (c) follows from (b) and Remark
\ref{adjoints}(iii). \end{proof}

In $\pp^3$ we get back the well known speciality bound $e\leq \frac
{(d-4)}2$, attained for instance by complete intersections on quadrics. A
complete classification of the curves in $\pp^3$ achieving this speciality
bound will be given in Proposition \ref {alpha minimal} (see also Remark
\ref{quadric1'}). \bigskip

\begin{rem} \rm Concerning the relations between the values $\partial
h_N(\alpha-1) = \alpha$ and $\partial h_N(\alpha)$ we can say the following:
\begin{itemize}
\item[(a)] by Theorem \ref{pr} we have $0 \le \alpha - \partial h_N(\alpha)
\le \alpha - (d-2-\alpha)(n-2)$;
\item[(b)] $\partial h_N(\alpha) = \alpha$ if $C$ is subcanonical, but not
conversely.
Indeed $C$ is subcanonical if and only if $\alpha d = 2\delta$. Moreover in
this case there is a double structure on $N$ which is a complete
intersection of type $(\alpha, \delta)$ (see
\cite {CGN});
\item[(c)] if $\alpha$ is minimal with respect to $d$ and $r$, namely
$\alpha = \lceil \frac{(d-2)(r-2)}{r-1} \rceil$
(see Corollary \ref{bound for alpha}(b)), then
$\alpha - \partial h_N(\alpha) \le n-2$, with equality if and only if
$\alpha = \frac{(d-2)(r-2)}{r-1}$. This is an easy consequence of Corollary
\ref{bound for alpha} (details to
the reader).
\end{itemize}
\end{rem}

We end this section by showing that when $C$ is arithmetically
Cohen-Macaulay (aCM) $\partial h_N$ can be computed from the Hilbert
function of a general hyperplane section of $C$ (and
conversely).

We recall first an easy fact.

\begin{rem}\label{nodes and hyperplane section}\rm Let $\Gamma$ be a
general hyperplane section of $C$. Then it is easy to see that
$$
h^0(\OO_{C}(j)) - h^0(\OO_{C}(j-1)) \ge h_\Gamma (j) \quad {\rm for \;
every} \quad j \in \ZZ
$$
with equality whenever $h^1(\mathcal I_{C}(j)) = 0$:
moreover equality holds for all $j \in \ZZ$ if and only if $C$ is aCM.

Hence by Remark \ref{adjoints}(iii) we have
$$
\partial h_N(j) \ge j+1-d+ h_\Gamma (d-2-j) \quad {\rm for \; every} \quad
j \le d-1
$$
with equality whenever $h^1(\mathcal I_{C}(j)) = 0$; moreover equality holds
for every $j \le d-1$ if and only if $C$ is aCM.
\end{rem}

\begin{ex}\label{aCM} {\rm  Assume that $C$ is aCM and let $\Gamma$ be a
general hyperplane section of $C$.
Then by Remark \ref{nodes and hyperplane section} we have
$$
\partial h_N(j) = j+1-d+ h_\Gamma (d-2-j) \quad {\rm for \; every} \quad j
\le d-1.
$$
It is also easy to see that the index of speciality of $C$ is
$$
e = {\rm max}\{j \in \ZZ \mid h^1(\II_\Gamma (j)) \not = 0\} - 1
$$
whence
$$
\alpha = d - 2 - {\rm max}\{j \in \ZZ \mid h^1(\II_\Gamma(j))\not = 0\}.
$$
Then, since $\partial h_N(j) = j+1$ for $0 \le j \le \alpha-1$ it is
immediate to write down
$\partial h_N(j)$ for every $j \in \ZZ$, in terms of $d$ and $h_\Gamma$.

For example if $r = 3$ and $C$ is a complete intersection of type $(s,
s+u)$, with $u \ge 0$, we have
$$
\partial h_\Gamma(j) = \begin{cases}
0&\text{ if } \; j<0 \\
j+1 &\text{ if } \;  0\leq j \le s-1\\
s &\text{ if } \;  s\leq j \le s+u-1\\
2s+u-1-j &\text{ if } \;  s+u\leq j \le 2s+u-2\\
0&\text{ if } \;  j\ge 2s+u-1
\end{cases}
$$
It follows (as expected):
$$
e = 2s+u-4
$$
whence
$$
\alpha = d+1-(2s+u)\quad
$$
and, with some calculations,

$$\partial h_N(j)=\begin{cases}
0&\text{ if } \;  j<0 \\
j+1 &\text{ if } \;  0\leq j< \alpha\\
\alpha-\binom{j-\alpha+1}2 &\text{ if } \; \alpha\leq j\leq d-s-u-2 =
\alpha+s-3\\
\binom s2+(d-s-1-j)(s-1) &\text{ if } \;  d-s-u-1\leq j\leq d-s-2\\
\binom{d-1-j}2 &\text{ if } \;  d-s-1\leq j\leq d-3\\
0&\text{ if } \; j>d-3
\end{cases}
$$
 }
 \end{ex}
\bigskip

\begin{rem} \label {spectrum} \rm The Hilbert function of $N$ is strictly
related with the {\it spectrum} $\ell_C$  of $C$, defined as (see \cite{SCH}):
$$
\ell_C (j) := \partial^2 h^0(\OO_C(j)) \quad {\rm for \; all}\quad j \in \ZZ
$$

Indeed from Remark \ref{adjoints}(iii) it follows easily  that
\begin{equation}\label{spectrum_and_hN}
\partial^2 h_N(t) = 1-\ell_{C}(d-1-t) \quad {\rm for} \quad \alpha \le t
\le d-1
\end{equation}

For example if $C$ is aCM and $\Gamma$ is a general hyperplane section of
$C$ the spectrum of $C$ is
$$
\ell_{C}(j) = \partial h_\Gamma (j) \quad {\rm for} \quad j \in \ZZ.
$$
as follows easily by Example \ref{aCM}.
\end{rem}
\bigskip

\section{Postulation of nodes and Castelnuovo's curves}\label{Nodes and
Castelnuovo's curves}

In this section we use Theorem \ref {pr} to produce a lower bound for
$\partial h_N$ in
terms of $d$ and $r$ and we show that this bound is sharp for curves of
maximal genus, that is
Castelnuovo's curves. Moreover we give some lower and upper bounds for the
length $\delta$ of $N$ depending on $d$, $r$ and $\alpha$, hence lower and
upper bounds for $g$ depending on
$d$, $r$ and $e$.

\begin{prop} \label{Castelnuovo curves} Set $\beta := \lceil
\frac{(d-2)(r-2)}{r-1}\rceil$ and
$q:= d-2-\beta$ and let $\phi_\beta$ be the function defined by:
$$
\phi_\beta(j) = \begin{cases}
0&\text{ if } \; j<0 \\
j+1 &\text{ if } \;  0\leq j \le \beta-1\\
(d-2-j)(r-2) &\text{ if } \;   \beta \le j \le d-3\\
0&\text{ if } \;  j\ge d-2
\end{cases}
$$
Then we have:
\begin{itemize}
\item[(a)] $\partial h_N (j) \ge \phi_\beta(j)$ for every $j \in \ZZ$;
\item[(b)] $\delta \ge \frac 12 (\beta(\beta+1)+(r-2)(q^2+q))$;
\end{itemize}

Moreover the following conditions are equivalent:
\begin{itemize}
\item [(i)] $\partial h_N (j) = \phi_\beta(j)$ for every $j \in \ZZ$;
\item [(ii)] $\delta = \frac 12 (\beta(\beta+1)+(r-2)(q^2+q))$;
\item [(iii)] $C$ is a Castelnuovo curve (that is a curve of maximal genus
among the non degenerate curves
of degree $d$ in $\pp^r$).
\end{itemize}
\end{prop}

\begin{proof}

By Corollary \ref {bound for alpha} we have $\beta \le \alpha$  and (a)
follows by Theorem \ref{pr}. Moreover (a) implies
$\delta = \sum_{j=0}^{d-3} \partial h_N(j) \ge \sum_{j=0}^{d-3} \phi_\beta
(j) = \frac 12 (\beta(\beta+1)+(r-2)(q^2+q))$,
whence (b).

The equivalence of (i) and (ii) follows immediately from (a) and (b).

To conclude the proof we recall that Castelnuovo's curves are characterized
by the equality
\begin{equation}\label{Castelnuovo}
g = \binom m 2 (r-1) + m\epsilon
\end{equation}
where
$$
d-1 = (r-1)m + \epsilon, \quad 0 \le \epsilon \le r-2.
$$
(see e.g. \cite {EH} or  \cite {C} for details).

Now equality in (b) is equivalent to
\begin{equation}\label{delta minimal}
g = \frac {(d-1)(d-2)}{2} - \frac 12 (\beta(\beta+1)+(r-2)(q^2+q)).
\end{equation}

It is easy to see that $q = d-2-\beta = [\frac {d-2}{r-1}]$, that is
$$
d-2 = (r-1)q + \eta, \quad 0 \le \eta \le r-2.
$$

This implies, by a straightforward calculation, that  (\ref {delta
minimal}) is equivalent to
\begin{equation}\label{Castelnuovo1}
g = \binom q 2 (r-1) + q(\eta + 1).
\end{equation}

Now if $\eta \le r-3$ we have $(m,\epsilon) = (q, \eta +1)$ and if $\eta =
r-2$ we have $(m,\epsilon) = (q+1, 0)$.
It follows that in both cases  (\ref {Castelnuovo1}) is equivalent to
(\ref {Castelnuovo}).
\end{proof}

\bigskip

\begin{rem} \rm Proposition \ref{Castelnuovo curves} shows that the lower
bound for $\partial h_N$ provided by Theorem \ref{pr} and the lower bound
for $\alpha$ given by Corollary \ref{bound for alpha} are sharp, since they
are attained by Castelnuovo's curves. However we will see that these bounds
can be achieved also by other curves in $\pp^3$, where the situation is
quite clear (see Section \ref{low degree}). In higher-dimensional spaces
several things are still unclear. For example even considering complete
intersections of quadrics we have the following situation.

Let $C$ be a complete intersection of $3$ quadrics in $\pp^4$. Then $d=8$
and $g=5$, so that $C$ is a Castelnuovo curve. Indeed a direct calculation
easily shows that
$$
\partial h_N :...,0,1,2,3,4,4,2,0,...
$$
coherently with Proposition \ref{Castelnuovo curves}.

On the other hand, the reader can easily realize (e.g. by using Example
\ref{aCM}) that complete intersections of quadrics in higher dimensional
projective spaces, which are not Castenuovo's curves, do not have minimal
descent of $r-2$ at any step after $\alpha+1$.
\end{rem}
\medskip

Also in $\pp^4$ more general complete intersections do not have the minimal
descent, as the following example shows.

\begin{ex} \rm Let $C$ be a complete intersection of type $(2,2,3)$ in
$\pp^4$.
Then $d=12$ and $g=13$, whence $\delta = 55-13=42$.

Since $\omega_{C}=\OO_{C}(2)$, by Remark \ref{adjoints} we have
$\alpha=d-3-2=7$.Then $\partial h_N(j)=j+1$ for $j<7$ while $\partial
h_N(j)=0$ for $j>d-3=9$; furthermore the sum of $\partial
h_N$ is $42$.

It remains to compute $\partial h_N(j)$ for $j=7,8,9$. Theorem \ref{pr}
says that $\partial h_N(8)\leq \partial h_N(7)-2$, $\partial h_N(9)\leq
\partial h_N(8)-2$, while $ \partial h_N(9)=2$ by proposition \ref{basic}.
Summing up, we see that the unique possibility for $\partial h_N$ is
 $$\partial h_N: ...,0,1,2,3,4,5,6,7,7,5,2,0,...$$

Observe that there is an intermediate step in which the descent is by $3$,
i.e. more than the minimum allowed $r-2=2$, whence, in particular, $C$ is
not a Castelnuovo curve. However $\alpha$
reaches the minimal value established by Corollary \ref {alpha minimal}.
\end{ex}
\medskip

\begin{prop} \label{bounds for delta}
Assume $g > 0$, set
$$
\begin{array}{l} p : = d-2-\alpha,\\
 \mu := \alpha - p(r-2),\\
\end{array}
$$
and write $\alpha$ as
$$
\alpha = (r-2)m+\nu, \quad 0 \le \nu \le r-3.
$$

Then $p > 0$, $\mu \ge 0$ and $\delta$ satisfies the following inequalities:
\smallskip

\begin{itemize}
\item[(a)] $\frac 12 (\alpha(\alpha+1)+(r-2)(p^2+p))\le \delta\leq \frac 12
(\alpha(\alpha+1)+(r-2)(p^2+p))+\mu(p-1)$;
\medskip

\item[(b)] $\delta\leq \frac 12 (\alpha(\alpha+1)+(r-2)(m^2+m))+\nu(m-1)$.
\end{itemize}
\end{prop}

\begin{proof} We have $p > 0$ by Remark \ref{adjoints}(v) and $\mu \ge 0$
by Corollary \ref {bound for alpha}.
Consider now the numerical functions $\phi$ and $\Phi$ defined by:
$$
\phi(j) = \begin{cases}
0&\text{ if } \; j<0 \\
j+1 &\text{ if } \;  0\leq j \le \alpha-1\\
(d-2 - j)(r-2) &\text{ if } \;   \alpha \le j \le d-3\\
0&\text{ if } \;  j\ge d-2
\end{cases}
$$
$$
\Phi(j) = \begin{cases}
\phi(j)&\text{ if } \; j\le \alpha-1\\
\phi(j) + \mu &\text{ if } \;   \alpha \le j \le d-4\\
\phi(j)&\text{ if } \;  j \ge d-3
\end{cases}
$$

By Theorem \ref{pr} and Proposition \ref{basic} we have that $\phi(j) \le
\partial h_N(j) \le \Phi(j)$ for every $j$. Moreover $\delta =
\sum_{j=0}^{d-3} \partial h_N(j)$, whence $\sum_{j=0}^{d-3} \phi(j) \le
\delta \le \sum_{j=0}^{d-3} \Phi(j)$. Then (a) follows by a straightforward
calculation.

To show (b) it is sufficient to show that $\frac 12 (r-2)(p^2+p)+\mu(p-1)
\le \frac 12 (r-2)(m^2+m)+\nu(m-1)$. Since $\mu \ge 0$ it is easy to see
that $m \ge p$ and $\mu = (r-2)(m-p)+\nu$. The conclusion follows by a
direct computation.
\end{proof}
\medskip

\begin{rem} \rm The bounds of Proposition \ref{bounds for delta} are
equivalent to bounds for the genus in terms of $d$, $r$ and $e$. The easy
transformation is left to the reader.
\end{rem}
\bigskip

\section{Postulation of nodes and flag conditions}

The sharpness of the previous results is rather flabby because as soon as
one takes curves which do not lie on surfaces of minimal degree, then
stronger restrictions seem to apply
(consider for example complete intersections of surfaces of degree at least
$3$ in $\pp^3$, see \ref{CI}).\par

So one is led to ask whether  more information about the flags of
subvarieties of $\pp^r$ containing $C$ could give more restrictive
constraints to the behavior of $\partial h_N$.

To obtain sharper information on the shape of $\partial h_N$
one may start with the formula
\begin{equation}\label{DeltahN1}
\quad\partial h_N(j)=j+1-h^0(\omega_{C}(j-d+3))+h^0(\omega_{C}(j-d+2))
\end{equation}
(see Remark \ref{adjoints}(iv)). This formula can be read as
$$
\partial h_N(j)=j+1-h^0(\OO_C(L+jD))+h^0(\OO_C(L+(j+1)D))
$$
where $D$ is a hyperplane divisor and $L: = K - (d-3)D$ is a divisor
corresponding to the invertible sheaf $\omega_{C}(-d+3)$. Then one is led
to study the growth of the linear series $|L+jD|$ with respect to $j$, that
is the growth {\it Hilbert function $h_L$ of $L$}, defined as:
$$
h_L(j) := h^0(\OO_C(L+jD)) \quad {\rm for \  every}\quad j \in \ZZ,
$$
which was introduced in  \cite {AS} and studied intensively by many authors.

With this notation (\ref{DeltahN1}) can be rewritten as
\begin{equation}\label{DeltahN1'}
\quad\partial h_N(j)=j+1-\partial h_L(j+1), \quad \OO_C(L) = \omega_{C}(-d+3)
\end{equation}
and any result on the step-by step growth of $h_L$ can be translated into
results
on the descent of $\partial h_N$ from $\alpha$ to $0$, and conversely.\par

Equivalently, one may start from equation
\begin{equation}\label{DeltahN2}
\quad \partial h_N(j)=j+1-d-h^0(\OO_{C}(d-j-3))+h^0(\OO_{C}(d-j-2))
\end{equation}
in Remark \ref{adjoints}(iii), which can be rewritten as
\begin{equation}\label{DeltahN2'}
\partial h_N(j)=j+1-d+\partial h_L(-j+1), \quad L : = (d-3)D
\end{equation}
Here information on the growth of $h_L$ are translated, reversing the
orientation, into information on the decrease of $\partial h_N$, starting
from $d-3$ and going backward to $\alpha$.
\smallskip

We will take both points of view in order to relate geometric properties of
the embedding $C\subset\pp^r$ to arithmetic properties of the functions
$h_L$ above.\par

\begin{rem} \rm In order to show how this method works we give an
alternative short proof (in characteristic zero) of Theorem \ref {pr}.
\smallskip

Let $M $ be a divisor corresponding to $\omega_{C}(j-d+3)$ and let $D$ be a
hyperplane divisor of $C$. Then by assumption one has $h^0(\OO_{C}(M-D))>0$
for all $j\geq \alpha+1$. Let $\Gamma$ be a general hyperplane section of
$C$ and look at the exact sequences:
$$ 0\to H^0 (\OO_{C}(M-D))\to H^0 (\OO_{C}(M))\to H^0(\OO_\Gamma) $$
$$0\to H^0 (\OO_{C}(M))\to H^0 (\OO_{C}(M+D))\to H^0(\OO_\Gamma) $$
Let $v=h^0(\OO_{C}(M))-h^0(\OO_{C}(M-D))$ and
$v'=h^0(\OO_{C}(M+D))-h^0(\OO_{C}(M))$ be the dimensions of the images in
$H^0 (\OO_\Gamma)$ of the two left hand side maps above.
Then $\Gamma$ is in uniform position with respect to the linear series
$|M|$ and $|D|$, so one may apply Castelnuovo's lemma of \S 1 in
\cite{CCD1}. It turns outthat $v'\geq v+r-1$. A direct computation shows
that the same inequality holds also when $j=\alpha$.

Now by equality (\ref{DeltahN1}) above, we obtain:
$$
\begin{array}{rcl}\partial h_N(j+1)-\partial h_N(j) &=&
h^0(\omega_{C}(j-d+4))-2h^0(\omega_{C}(j-d+3))+\\
& &+h^0(\omega_{C}(j-d+2))-1\\
&\geq& v'-v-1
\end{array}
$$
and the claim follows.
\end{rem}\medskip

Now we go back to our main problem, namely to take into account some flag
conditions in order to get sharper information on $\partial h_N$. As far as
we know, such a theory is available only for curves in $\pp^3$, where the
knowledge of a number $s$ such that $C$ is contained in no surfaces of
degree smaller than $s$ influences the growth of the dimensions of linear
series on the curve, under the addition of multiples of $D$. Indeed let us
recall a corollary of the main
technical result of \cite{CCD2}, appearing in \cite{CC} (see also
\cite{SCH} for generalizations to non-integral curves):
\smallskip

\begin{thm} \label{cc} Let $C\subseteq \pp^3$ be a smooth connected curve
and let $L$ be a divisor on $C$. Let $t$ and $j$  be integers such that:
$$h^0(\II_C(t-1))=0,\leqno{(a)}$$
$$\binom{t-1}2 < \partial h_L(j) < d -\binom t2. \leqno{(b)}$$
Then:
$$\partial h_L(j+1) \ge \partial h_L(j) + t.$$
\end{thm}

\proof See \cite{CC}, Theorem 2.3.\qed\medskip

Notice that the assumption of the above theorem implies $\binom{t-1}2 < d
-\binom t2$, hence also $d>(t-1)^2+1$, whence by  Laudal's Lemma also the
general hyperplane section of $C$ is not
contained in any curve of degree $t-1$.\par
\smallskip

We point out the consequences of the previous theorem for the Hilbert
function of $N$, in the following:

\begin{thm}\label{p3} Let $C\subseteq \pp^3$ be a smooth connected curve
not lying on any surface of degree smaller than $s$. Then  for all $j$ with
$\alpha\leq j\leq d-2$, if $t$ is any integer satisfying:
$$
t\leq s \leqno{(a)}
$$
$$
\binom{t-1}2 < j+1-\partial h_N(j)< d -\binom t2\leqno{(b)}
$$
one has:
$$
\partial h_N(j+1)\leq \partial h_N(j)-t+1.
$$
\end{thm}

\proof The claim follows easily by (\ref{DeltahN1'}) and Theorem \ref{cc}.
\qed\medskip

We want to explain roughly how one can use Theorem \ref{p3} to understand
the behavior of $\partial h_N(j)$, for $j\in [\alpha+1,d-3]$. Indeed
$\alpha$ is the maximum for the function
$\partial h_N$, achieved for $j=\alpha-1$; in the next few steps, when $j$
is a little bit bigger than $\alpha$, then the difference $j+1-\partial
h_N(j)$ can be small and one may apply theorem \ref{p3} only for small (but
nevertheless increasing) $t$; so the rate of the descent
of $\partial h_N$ is small, but increasing. Eventually $j+1-\partial
h_N(j)$ becomes bigger than  $\binom {s-1}2$, and one is allowed to apply
the theorem taking $s=t$ (of course provided
that $j+1-\partial h_N(j)$ is still smaller than $d-\binom s2$): here
$\partial h_N(j)$ is forced to decrease by $s-1$ (at least) at any steps.
Finally when $j$ approaches $d-2$, then $j+1-\partial h_N(j)$ becomes
bigger than $d-\binom s2$ and once again one must take $t<s$, in order to
apply the theorem: in other words the rate of the descent of $\partial h_N$
is
allowed to become smoother.\par

The following picture synthesizes the previous discussion on the shape of
$\partial h_N$:

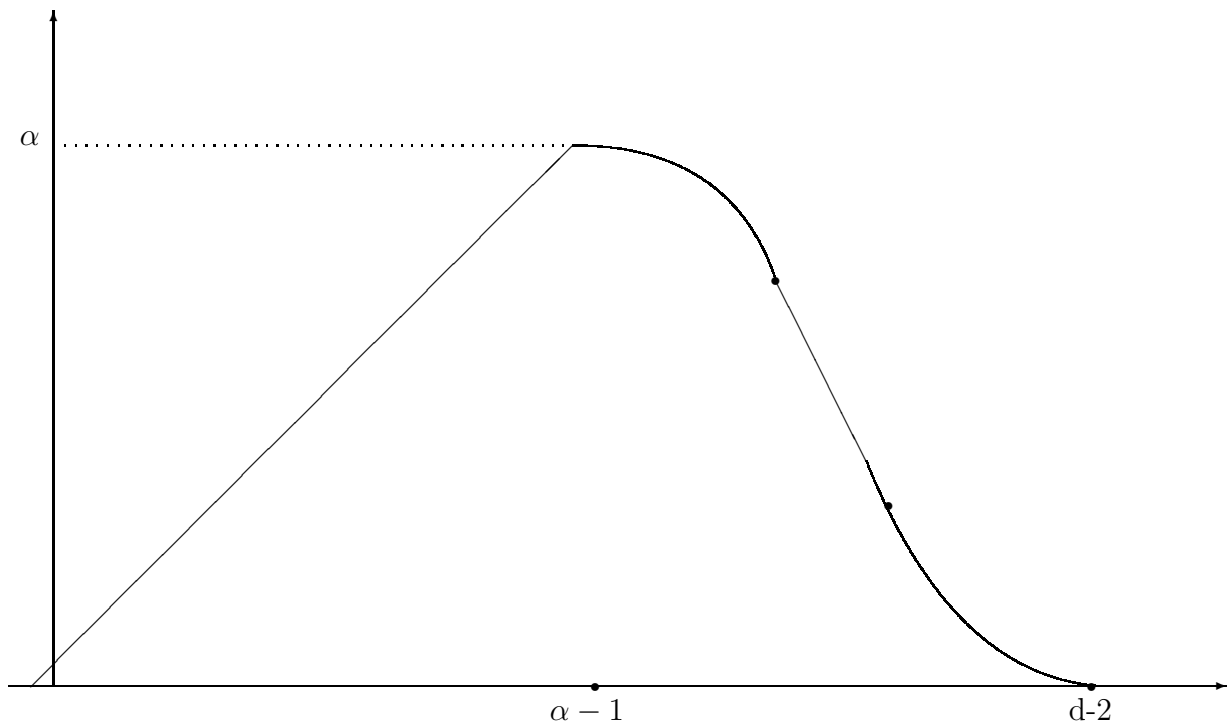
\begin{picture}(260,300)(-1,-10)
\thinlines \setlength{\unitlength}{3mm}

\put(-2,0){\vector(1,0){54}} \put(0,0){\vector(0,1){30}}



\put(22,-1.5){$\alpha-1$}
\put(24,0){\circle*{.3}}
\put(-1,0){\line (1,1){24}}
\put(0,0){\qbezier[1000](23,24)(30,24)(32,18)}
\put(32,18){\circle*{.3}}
\put(32,18){\line (1,-2){4}}
\put(0,0){\qbezier[1000](36,10)(40,0)(47,0)}
\put(37,8){\circle*{.3}}
\put(45,-1.5){d-2}
\put(46,0){\circle*{.3}}
\put(0,0){\qbezier[50](0,24)(12,24)(23,24)}
\put(-1.5,24){$\alpha$}
\end{picture}

\begin{ex}\label{IC46A} \rm Let us see the minimal descent of $\partial
h_N$ in the following numerical case: $d=24$, $s=4$, $g=73$, $\alpha=15$
(here $\delta=180$).

The function $\partial h_N$ increases by $1$ at any step from $0$ to
$\alpha-1=14$. Then:\par

- for $j=15$, we get $15=\alpha\geq \partial h_N(15)$. Assume the maximum
is attained; then for $t=2$, one sees that the inequalities (a),(b) of the
previous theorem are satisfied (in fact $t=2$ is the maximum which
satisfies the left hand side) so $\partial h_N$ is forced to decrease at
least by $1$;\par

- for $j=16$, we get then $14\geq \partial h_N(16)$. Assume the maximum is
attained; then for $t=3$, one sees that the inequalities (a),(b) of the
previous theorem are satisfied (in fact $t=3$ is the maximum which
satisfies the left hand side) so $\partial h_N$ is forced to decrease  at
least by $2$;\par

 - for $j=17$, we get then $12\geq \partial h_N(17)$. Assume the maximum is
attained; then for $t=4$, one sees that the inequalities (a),(b) of the
previous theorem are
satisfied (in fact $t=4$ is the maximum which satisfies the left hand side)
so $\partial h_N$ is forced to decrease  at least by $3$;\par

 - for $j=18$, we get then $9\geq \partial h_N(18)$. Assume the maximum is
attained; then  the inequality (b) of the previous theorem holds for $t=5$,
but $5>s$; we may just conclude that
$\partial h_N$ is forced to decrease  at least by $3$;\par

 - for $j=19$, we get then $6\geq \partial h_N(19)$. Assume the maximum is
attained; then  for $t=4$, one sees that the inequalities (a),(b) of the
previous theorem are satisfied (in fact, now $t=4$ is the maximum which
satisfies the {\it right} hand side) so $\partial h_N$ is forced to
decrease  at least by $3$;\par

 - for $j=20$, we get then $3\geq \partial h_N(20)$. Assume the maximum is
attained; then  for $t=3$, one sees that the inequalities (a), (b) of the
previous theorem are satisfied (in
fact $t=3$ is the maximum which satisfies the right hand side) so $\partial
h_N$ is forced to decrease  at least by $2$;\par

 - for $j=21$, we get then $1\geq \partial h_N(21)$. Here the maximum {\it
must be} attained, for $21=d-3$ and we know that $\partial h_N(d-3)\ge 1$
for all curves in $\pp^3$.\par

Finally, if the maximum is attained in all steps, the function $\partial
h_N$ reads:
$$
 ...,0,1,2,3,4,5,6,7,8,9,10,11,12,13,14,15,15,14,12,9,6,3,1,0,...
$$
and the sum of this function is exactly $180=\delta$, as it should be for
the Hilbert function of nodes of a concrete curve.
\end{ex}

A similar conclusion holds taking the reverse point of view and using (\ref
{DeltahN2'}) instead of (\ref {DeltahN1'}).

\begin{thm}\label{p3bis} Let $C\subseteq \pp^3$ be a smooth connected curve
not lying
on any surface of degree smaller than $s$. Then for all $j$ with
$\alpha\leq j\leq d-2$, if $t$
is any integer satisfying:
$$
t\leq s \leqno{(a)}
$$
$$
\binom{t-1}2 < d-j-1+\partial h_N(j)< d -\binom t2 \leqno{(b)}
$$
one has:
$$
\partial h_N(j-1)\geq \partial h_N(j)+t-1.
$$
\end{thm}

\proof The claim follows from (\ref{DeltahN2'}) and Theorem \ref {cc}.
\qed\medskip

\begin{ex}\label{IC46B} \rm Let us take the previous numerical example,
starting this time from $j = d-3$ and going backwards; assume still $d=24$,
$s=4$, $g=73$, $\alpha=15$ (here $\delta=180$).

The function $\partial h_N$ increases by $1$ at any step from $0$ to
$\alpha-1=14$, so that also $\partial h_N(14)=15$.\par

- for $j=21=d-3$, we must have $\partial h_N(21)\ge 1$. Assume that the
minimum is attained;
then for $t=3$, one sees that the inequalities (a),(b) of Theorem
\ref{p3bis} are satisfied (in fact $t=3$ is the maximum which satisfies the
left hand side) so $\partial h_N$ is forced to increase  at least by
$t-1=2$, going one step back;\par

- for $j=20$, we get then $\partial h_N(20)\geq 3$. Assume the minimum is
attained; then for $t=4$, one sees that the inequalities (a),(b) of the
previous theorem are satisfied (in fact $t=4$ is the maximum which
satisfies the left hand side) so $\partial h_N$ is forced to increase  at
least by $3$, going one step back;\par

- for $j=19$, we get then $\partial h_N(19)\geq 6$. Assume the minimum is
attained; then for $t=5$, one sees that condition (b) of the previous
theorem is satisfied, but condition (a)
fails. All we can say, taking again $t=4$, is that $\partial h_N$ is forced to
increase  at least by $3$, going one step back;\par

- for $j=18$, we get then $\partial h_N(18)\geq 9$. Assume the minimum is
attained; then for $t=5$, one sees that condition (b) of the previous
theorem is satisfied, but condition (a) fails. All we can say, taking again
$t=4$, is that $\partial h_N$ is forced to increase  at least by $3$, going
one step back;\par

 - for $j=17$, we get then $\partial h_N(18)\geq 12$. Assume the minimum is
attained; then  for $t=3$, one sees that the inequalities (a),(b) of the
previous theorem are satisfied (in fact, now $t=3$ is the maximum which
satisfies the {\it right} hand side) so $\partial h_N$ is forced
to increase  at least by $2$, going one step back;\par

 - for $j=16$, we get then $\partial h_N(16)\geq 14$. Assume the minimum is
attained; then  for $t=2$, one sees that the inequalities (a),(b) of the
previous theorem are
satisfied (in fact, now $t=3$ is the maximum which satisfies the right hand
side) so $\partial h_N$ is forced to increase  at least by $1$, going one
step back;\par

 - for $j=15=\alpha$, finally we get then $\partial h_N(21)\geq 15=\alpha$.
Here the minimum {\it must be} attained.\par

Finally, if the minimum is attained in all steps, the function $\partial
h_N$ reads:
$$
 ...,0,1,2,3,4,5,6,7,8,9,10,11,12,13,14,15,15,14,12,9,6,3,1,0,...
$$
which are exactly the same values as before.
\end{ex}

Notice that the numbers $d=24$, $s=4$, $g=73$, $\alpha=15$ above are not
random:
they are the numbers of a complete intersection $C$ of type $(4,6)$, see
Example
\ref{aCM}.\par
\bigskip

\section{Minimal node functions and Halphen curves}

In this section  we use the previous results to construct a new function which
bounds $\partial h_N$ from below, and we show that this bound is sharp for the
so-called Halphen's curves (see Theorem \ref{Halphen}).
Further examples will be given in the next section.

\begin{defn} \label{mnf} \rm (a) Given integers $d,s,n$ with $d \ge 3$,
$s\geq 2$ and $n\geq 3$ we define a numerical function
$\overline \Delta$ in the interval $[0,d-3]$ by descending induction
(motivated by Theorems \ref {pr} and \ref {p3bis})
as follows:
\smallskip
(1)\ \ \ $\overline \Delta(d-3):=n-2$ \\ (2)\ \ \ for $d-4 \ge j\geq 0$,
$$
\overline \Delta (j-1) := \overline \Delta (j)+\max \{t-1, n-2\}
$$
where
$$
t := \max \left \{ i\leq s \left \vert  \binom{i-1}2 \right. <
d-j-1+\overline \Delta(j) < d -\binom i2\right \}
$$
(we agree that the maximum of the empty set is $-\infty$). \par
\medskip

(b) Given integers $d,s,n,\alpha$ with $d \ge 3$, $s\geq 2$, $n\geq 3$ and
$\alpha \ge 0$ we define the {\it minimal node function} $\Delta$ with
respect to $d,s,n,\alpha$ as follows:
$$
\Delta(j) :=
\begin{cases}
0 & {\rm if}  \quad j<0 \quad {\rm or}\quad j>d-3\\
j+1& {\rm if}  \quad 0\leq j \leq \alpha-1\\
\overline \Delta(j) & {\rm if} \quad \alpha \le j \le d-3
\end{cases}.
$$
We set
$$\delta(\Delta):=\sum_{j=0}^\infty \Delta(j)$$
\end{defn}
\medskip

The next theorem shows that minimal node functions provide lower bounds for
$\partial h_N$, which are sharp in some significant cases, as we shall see
later.\par
\medskip

\begin{thm}\label{arit} Let $C\subseteq \pp^3$ be a smooth connected curve
of degree $d$,  lying on no surface of degree smaller than $s$ ($s\ge 2$).
Let $\Delta$ be the minimal node function, with respect to $d$, $s$, $n$,
$\alpha$. Then $\partial h_N(j)\geq \Delta(j)$ for all $j \in \ZZ$.
\end{thm}

\begin {proof} Notice that $\partial h_N(j)=\Delta(j)$ for all $j$ outside
of the interval $[\alpha, d-4]$, and in particular $\partial
h_N(d-3)=\Delta(d-3)=n-2$. We prove the non-trivial part of our statement
by descending induction on $j$, starting with $j = d-3$; namely we show
that, for $\alpha+1\leq j\leq d-3$,  $\partial h_N(j)\geq \Delta(j)$
implies $\partial h_N(j-1)\geq \Delta(j-1)$.

Assume then that $\partial h_N(j)\geq \Delta(j)$ for some  $j$  with
$d-3\geq j\geq \alpha+1$. Set

$$
t: = \max \left \{i\leq s \left | \binom{i-1}2 < d-j-1+\Delta(j)< d -\binom
i2\right. \right \}
$$
and similarly

$$
t': = \max \left \{i\leq s \left | \binom{i-1}2 < d-j-1+\partial h_N(j)< d
-\binom i2 \right. \right \}.
$$
\smallskip

We observe that $t' \ge 1$. Indeed  $d-j-1+\partial h_N(j) < d -\binom 12 =
d$, because $\partial h_N(j) \le \alpha < j+1$, and $d-j-1+\partial h_N(j)
> \binom 02 = 0$ because $\partial h_N(j) > 0$.

If $t-1\le n-2$ we have $\Delta(j-1) = \Delta(j) +n-2$, and since $\partial
h_N(j-1)\geq\partial h_N(j)+n-2$ by Theorem \ref{pr}, the conclusion is clear.

If $t-1 > n-2$ we argue by contradiction, assuming that $\partial h_N(j-1)
< \Delta(j-1)$. Since $\Delta(j-1)=\Delta(j)+t-1$ by definition, while
$\partial h_N(j-1)\geq\partial h_N(j)+t'-1$ by Theorem \ref{p3bis}, the
above inequality and the induction assumption imply
\begin{equation}\label{not0}
t-t' > \partial h_N(j)-\Delta(j)\ge 0
\end{equation}
whence necessarily $t'+1\le t \le s$. Now the definition of $t'$ guarantees
that at least one of the following inequalities holds:
\begin{equation}\label{not1}
\binom {t'}2 \geq d-j-1 + \partial h_N(j)
\end{equation}
\begin{equation}\label{not2}
 d-j-1 + \partial h_N(j) \ge d - \binom {t'+1}2
\end{equation}
On the other hand by the definition of $t$, the inequality $t' \le t-1$ and
the induction hypothesis we have:
$$
\binom{t'}2 < d - j - 1 + \Delta(j) \le d - j - 1 + \partial h_N(j)
$$
whence (\ref{not1}) is false and (\ref{not2}) must be true.

Set $u : = t - t' - 1 \ge 0$. Then by (\ref{not2})  we have:
$$
\partial h_N(j) \ge j+1- \binom{t'+1}2
$$
while
$$
-\Delta (j) > -j-1+\binom t2
$$
by definition. Then by (\ref{not0}) we get
$$
u \ge  \partial h_N(j) - \Delta(j) > \binom t2-\binom{t'+1}2
= \frac {(2t'+1)u + u^2}2 \ge 0
$$
whence $u > 0$ and $1 > 2t' + u$. This is a contradiction because $t' > 0$.
\end{proof}
\medskip

In view of Theorem \ref{arit} it is natural to give the following Definition.

\begin{defn} \label{delta2}\rm
Let $C\subseteq \pp^3$ be a curve of degree $d$ not lying on any surfaces
of degree $s$ $(s\ge 2)$. We say that $C$ has {\it minimal nodal
postulation} when $\partial h_N$
coincides with the minimal node function with respect to $d, s, n, \alpha$.
\end{defn}

From the above theorem we get the following two Corollaries which give an
upper bound for the genus and a lower bound for $\alpha$ (whence an upper
bound for the speciality). These bounds are sharp in some important cases
(see Proposition \ref{HGP}, Theorem \ref{Halphen} and section \ref{low
degree}).

\begin{corol}\label{bound for g} Let the notation and the assumptions be as
in Theorem \ref{arit}.
Then:
\begin{itemize}
\item[(i)] $\delta \ge \delta (\Delta)$ and equality holds if and only if
$C$ has minimal node postulation;
\item[(ii)] $g \le \frac {(d-1)(d-2)}2 - \delta(\Delta)$ and equality holds
if and only if $C$ has minimal nodal postulation.
\end{itemize}
\end{corol}
\begin{proof} It is an immediate consequence of Theorem \ref {arit}.
\end{proof}

\begin{corol}\label{bound for alpha 2} Let the notation and the assumptions
be as in Theorem \ref{arit} and let $\overline \Delta$ be as in Definition
\ref{mnf}. Set
$$
\overline \alpha := \min \{j \ge 0\; \vert \; \overline \Delta(j) \le j\}.
$$
Then $\alpha \ge \overline \alpha$ and $e \le d-3-\overline \alpha$.
\end{corol}

\begin{proof} Let $\Delta$ be the minimal node function with respect to
$d$, $s$, $n$ and $\alpha$. Then by Theorem \ref{arit} and Proposition
\ref{basic} we have $\Delta(\alpha) \le \partial h_N(\alpha) \le \alpha$,
whence the conclusion by the definition of $\Delta$.
\end{proof}
\medskip

\begin{rem}\rm
(i) Let us observe that once we fix $d$, $s$ $n$ and $\alpha$, the minimal
nodal function $\Delta$ with respect to these values is completely
determined, so is its sum $\delta(\Delta)$.
Hence the upper bound for $g$ given by Corollary \ref {bound for g} depends
only on these numbers, or, equivalently, on $s$, $d$ $n$ and $e$.

Likewise if we fix $d$, $s$ and $n$ the lower bound $\bar \alpha$ for
$\alpha$ given by Corollary \ref {bound for alpha 2} is completely
determined by $d$ and $s$.
\smallskip

(ii) We will be able to compute explicitly $\bar \alpha$ under the
assumption $d > s(s-1)$ and $n = 3$ (see Proposition \ref {HGP}), and from
this we will be able to classify the curves with minimal nodal postulation
and minimal $\alpha$ in the given range: they are, of course,
the well-known Halphen curves (see Theorem \ref {Halphen}).

The study of curves with minimal nodal postulation but non-minimal $\alpha$
seems to be much harder; in particular it is not clear whether the above
bound for the genus is sharp. We will give some partial results for $d >
s(s-1)$ and  $s \le 4$ (see Section \ref {low degree}).

Our results and examples seem to indicate that the classification of curves
with minimal nodal postulation is not straightforward and could be a real
challenge.
\end{rem}

Before coming to the main results of this section we want to say something
on complete intersections. Observe that our previous examples \ref{IC46A}
and \ref{IC46B} show that
complete intersection curves of type $(4,6)$ have minimal nodal
postulation. More generally the following example shows that this is true
for all complete intersections.

\begin{ex}\label{CI} \rm Let $C\subseteq \pp^3$ be a complete intersection
of type $(s,s+u)$, with $s \ge 2$ and $u\geq 0$. Then $d=s^2+su$, $s$ has
the usual meaning, $\alpha=d-2s-u+1$ (see Example \ref {aCM}). Moreover
$\partial h_N(d-3) = 1$, whence $n=3$. The minimal nodal function $\Delta$
with respect to $d, \alpha, s$ as above and $n = 3$ can be computed
easily from the definition and is:
$$
\Delta(j)=\begin{cases}
0&\text{ if } \;  j<0 \\
j+1 &\text{ if } \;  0\leq j< \alpha\\
\alpha-\binom{j-\alpha+1}2 &\text{ if } \; \alpha\leq j\leq d-s-u-2 =
\alpha+s-3\\
\binom s2+(d-s-1-j)(s-1) &\text{ if } \;  d-s-u-1\leq j\leq d-s-2\\
\binom{d-1-j}2 &\text{ if } \;  d-s-1\leq j\leq d-3\\
0&\text{ if } \; j>d-3
\end{cases}
$$

By Example \ref{aCM} we have $\Delta = \partial h_N$, that is $C$ has
minimal nodal postulation.
\end {ex}

Now we want to show that minimal nodal postulation is achieved not only by
complete intersection
curves but also by Halphen's curves. We begin by showing a lower bound for
$\alpha$. This bound is well-known, but we give a new proof as a
consequence of our theory of minimal nodal functions.

\begin{prop}\label{HGP} Let $C \subseteq \pp^3$ be a curve of degree $d$
not lying on any surface of degree $< s$. Assume $s \ge 2$ and $d >
s(s-1)$. Let $(k,u)$ be the unique pair of integers such that $d = s^2+us -
k$, $0 \le k \le s-1$ and $u \ge 0$.
Assume further that $C$ is not a complete intersection (i.e. $k \ge 1$).
Then  we have:
\begin{itemize}
\item[(i)] $\bar\alpha \geq d-2s-u+2$, whence $\alpha \ge d-2s-u+2$ (where
$\bar\alpha$ is
defined in Corollary \ref {bound for alpha 2});
\item[(ii)] if $\alpha = d-2s-u+2$, then the minimal nodal function with
respect to $d$, $s$, $\alpha$ and $n=3$ is the following:
$$\Delta(j) =\begin{cases} 0&\text{ if } \;  j < 0 \\
j+1 &\text{ if } \;  0\leq j< \alpha\\
\alpha-\binom{j-\alpha+1}2 &\text{ if }\alpha\leq j\leq\alpha+k-2\\
\alpha-\binom{j-\alpha+1}2-(j+2-\alpha-k) &\text{ if }\alpha+k-1\leq
j\leq\alpha+s-4 = d-s-u-2\\
\binom s2+(d-s-1-j)(s-1) &\text{ if } \;  d-s-u-1\leq j\leq d-s-2\\
\binom{d-1-j}2 &\text{ if } \;  d-s-1\leq j\leq d-3\\
0&\text{ if } \; j>d-3
\end{cases}
$$
\end{itemize}
\end{prop}

\begin {proof}
To prove (i) we use Corollary \ref{bound for alpha 2}. A lengthy and
tedious but elementary calculation allows to compute the function
$\overline \Delta$ (see Definition \ref{mnf})
with respect to $d$, $s$ and $n = 3$. Going backwards step by step starting
with $d-3$ and using repeatedly the assumption $d > s(s-1)$ we get:
$$
\overline \Delta(j) =
\begin{cases}
\binom{d-1-j}2 &\text{ if } \;  d-s-1\leq j\leq d-3\\ \\
\binom s2+(d-s-1-j)(s-1) &
\begin{array}{ll}
\text{ if } &d-s-u-1\\
&\leq j\leq d-s-2
\end{array}\\ \\
\binom s2+u(s-1)+\displaystyle \sum_{t=1}^{d-s-u-1-j} (s-1-t) &
\begin{array}{ll}
\text{ if }& \;  d-2s-u+k+1\\
 &\le j\le d-s-u-2
\end{array}\\ \\
\binom s2+u(s-1)+\binom{s-1}2-\binom k2+\displaystyle
\sum_{t=2}^{d-2s-u+k+1-j} (k-t) &
\begin{array}{lcr}\text{ if }&  d-2s-u+2\le j &\\
 &\le d-2s-u+k& \end{array}\\ \\
\binom s2+u(s-1)+\binom{s-1}2-k+1 + d-2s-u+2 -j &\text{ if } \;  0 \le j
\le d-2s-u+1\\
\end{cases}
$$
If $k\ge2$ we have $\overline \Delta(d-2s-u+2) = d-2s-u+2$, and if $k=1$ we
have $\overline \Delta(d-2s-u+2) = d-2s-u+1$. Then (i) follows from
Corollary \ref{bound for alpha 2}.

The proof of (ii) is another tedious but elementary calculation using (i)
and the definition of $\Delta$.
\end {proof}

\begin{rem} \rm Since $e = d-3-\alpha$ it is immediate to see that the
bound of Proposition
\ref{bound for alpha 2} is equivalent to the Halphen-Gruson-Peskine bound
$e \le \frac ds+s-4$.
The first modern proof of this bound for $e$ is in \cite {GP}.

Another proof of the bound for $e$ can be found in \cite{SCH}, where the
Halphen-Gruson-Peskine bound is generalized to non-integral curves, using a
fine study of the spectrum.
\end{rem}

\begin{thm}\label{Halphen}
Let $C \subseteq \pp^3$ be a curve of degree $d$ not lying on any surface of
degree $<s$. Assume $s \ge 2$ and $d > s(s-1)$. Let $(k,u)$ be the unique
pair of
integers such that $d = s^2+us - k$, with $0 \le k \le s-1$ and $u \ge 0$. Then
the following are equivalent:
\begin{itemize}
\item[(i)] C lies on a surface $S$ of degree $s$ and it is either a complete
intersection of $S$ with a surface of degree $s+u$ (if $k = 0$), or it is
linked on $S$
to a plane curve of degree $k$, if $k > 0$;
\item[(ii)] if $k = 0$ the function $\partial h_N$ is as in Example \ref{CI}.
If $k > 0$ we have
$\alpha = d - 2s - u + 2$ and
$$
\partial h_N(j) =
\begin{cases}
0&\text{ if } \;  j < 0 \\
j+1 &\text{ if } \;  0\leq j< \alpha\\
\alpha-\binom{j-\alpha+1}2 &\text{ if }\alpha\leq j\leq\alpha+k-2\\
\alpha-\binom{j-\alpha+1}2-(j+2-\alpha-k) &\text{ if }\alpha+k-1\leq
j\leq\alpha+s-4 = d-s-u-2\\
\binom s2+(d-s-1-j)(s-1) &\text{ if } \;  d-s-u-1\leq j\leq d-s-2\\
\binom{d-1-j}2 &\text{ if } \;  d-s-1\leq j\leq d-3\\
0&\text{ if } \; j>d-3
\end{cases};
$$
\item[(iii)] $C$ has minimal nodal postulation with respect to $d,s$, $n =
3$ and
$$\alpha =
\begin{cases} d - 2s - u + 2&\text{ if } k > 0 \\
d - 2s - u + 1 &\text{ if } k = 0
\end{cases}
$$
In other words: if $\Delta$ is the minimal node function with respect to
$d$ $s$ $\alpha$ and $n = 3$ we have:
$$
\Delta(j) = \partial h_N(j) \quad {\rm for \; every} \; j \in \ZZ;
$$
\item[(iv)] $C$ has maximal genus in the set of all curves of degree $d$
not lying
on any surface of degree $< s$.
\end{itemize}
\end{thm}

\begin{proof} (i) $\Rightarrow$ (ii). If $k = 0$ follows immediately by
Example \ref{CI}.

Assume now $k > 0$ and let $Y$ be the plane curve linked to $C$ on $S$.
Then $\alpha=d-2s-u+2$, since $\omega_{C}$ is cut by surfaces of degree
$2s+u-4$ passing
through $Y$ (outside the intersection $Y\cap C$). Moreover $C$ is aCM, and
hence to
compute $\partial h_N$ we can use Example \ref {aCM}. Now the  general
hyperplane
section $\Gamma$ of $C$ is directly linked in a complete intersection
$(s,s+u)$ with the
general hyperplane section of $Y$, which is collinear of degree $k$. Then
by a well
known relation between the Hilbert functions of linked zero-dimensional
schemes (see
\cite {DGO}) we have
\begin{equation}\label{linked to collinear}
\partial h_\Gamma(j) =
\begin{cases}
0&\text{ if } \; j<0 \\
j+1 &\text{ if } \;  0\leq j \le s-1\\
s &\text{ if } \;  s-1\leq j \le s+u-1\\
2s+u-1-j &\text{ if } \;  s+u\leq j \le 2s+u-k-2\\
2s+u-2-j &\text{ if } \;  2s+u-k-1\leq j \le 2s+u-3\\
0&\text{ if } \;  j\ge 2s+u-2
\end{cases}
\end{equation}
On the other hand by Example \ref{aCM} we have
$$
\partial ^2 h_N(j) = 1 - \partial h_\Gamma(d-1-j)
$$
and a straightforward calculation shows (ii).

(ii) $\iff$ (iii). If $k = 0$ follows from Example \ref{CI}.
If $k>0$ use Proposition \ref{HGP}.

(iii) $\iff$ (iv). It is an easy consequence of Theorem \ref{arit} and
Proposition \ref{HGP}.

(ii) $\Rightarrow$ (i). Assume first that $k > 0$ and let $C'$ be a curve
of degree
$d$ linked to a plane curve of degree $k$ by a complete intersection of
type $(s,s+u)$.
Let $N'$ be the set of nodes of a general plane projection of $C'$. Then
$\partial h_N =
\partial h_{N'}$, whence by Remark \ref{adjoints}(iii) $h^0(\mathcal O_C(j)) =
h^0(\mathcal O_{C'}(j))$ for all $j \in \ZZ$. This easily implies, since
$C'$ is aCM,
that $h^0(\mathcal I_C(j)) - h^0(\mathcal I_{C'}(j)) = h^1(\mathcal
I_C(j))$, whence, in
particular, $h^0(\mathcal I_C(j)) \ge h^0(\mathcal I_{C'}(j))$. By easy degree
considerations this implies that $h^0(\mathcal I_C(j)) =  h^0(\mathcal
I_{C'}(j))$ for
$j \le s+u-1$, whence $h^1(\mathcal I_C(j)) = 0$ for $j \le s+u-1$. In
particular $C$ is
contained in an irreducible surface $S$ of degree $s$, and a surface $T$ of
degree $s+u$
not containing $S$.

Let now $\Gamma$ and $\Gamma'$ be the general hyperplane sections of $C$
and $C'$
respectively.
Then by Remark \ref {nodes and hyperplane section} it follows that
$h_{\Gamma'}(j) \ge h_\Gamma (j)$ for every $j$. We want to show that
$h_{\Gamma'} =
h_\Gamma$. Assume this is not true and set $t := \min \{j \in \ZZ \vert
h_{\Gamma'}(j) >
h_\Gamma (j)\}$. By the above considerations we have that $t \ge s+u$ and
that $\partial
h_\Gamma$ is strictly decreasing for $j \ge t$ until it vanishes. Now $\partial
h_{\Gamma'}$ is explicitly given in (\ref{linked to collinear}) (with
$\Gamma$ replaced
by $\Gamma'$), and since $\sum  \partial h_{\Gamma'}(j) = d = \sum \partial
h_\Gamma(j)$
an easy calculation leads to a contradiction.

Then $h_{\Gamma'} = h_\Gamma$, whence $C$ is aCM by Remark
\ref{nodes and hyperplane section} and $\Gamma$ is linked to a collinear
scheme in a
complete intersection of type $(s,s+u)$. Then it is easy to see that $S
\cap T = C \cup
C''$, where $C''$ is a plane curve of degree $k$.

If $k=0$ the proof is similar to (and easier than) the previous one and is
left to
the reader.
\end{proof}
\medskip

\begin{ex}\label{Cast} \rm There are some (few) examples in which the
minimal node function is the effective Hilbert function of the nodes of a
curve coming from $\pp^r$, $r\geq 4$. This is the case, for instance, of a
curve $C\subseteq \pp^3$ which is a general projection of a Castelnuovo
curve in $\pp^r$, at least when $d> (r-2)(r-1)$.\par\noindent

Indeed every Castelnuovo  curve $C' \subseteq \pp^r$ with $d > 2r$  lies on
a surface
of (minimal) degree $r-1$ (see \cite  {C} or \cite {EH}), which is also the
minimal
degree of a surface containing the projection $C$ of $C'$ to $\pp^3$, if
$d> (r-2)(r-1)$.
Then in the definition of minimal node function, we get here $t\leq s=r-1$
so that
$\max\{t-1, r-2\}=r-2$ hence, by construction (\ref{mnf} (2)), the minimal
node function
satisfies:
$$
\Delta(j-1) = \Delta(j)+r-2 \quad {\rm for} \quad \alpha \le j \le d-3.
$$
Now just observe that Castelnuovo's curves reach exactly this bound by
Proposition \ref{Castelnuovo curves}.

In order to get more refined bounds, similar to the ones in $\pp^3$,
it seems reasonable to take into account some flag conditions, such as, for
example, the
ones used in \cite {CCD2} to bound the genus.
\end {ex}

\begin{rem} \rm The curves of the previous example are the first examples of
non-aCM curves attaining the lower bound stated in Theorem \ref{arit}.
Indeed they  are
not even linearly normal.

Observe also that for these curves we have
$\alpha = \lceil \frac{(d-2)(r-2)}{r-1} \rceil$, whence $\alpha$ is not
minimal with
respect to $d$ and $s= r-1 >2$.

We will see in the next section that if $\alpha$ is not minimal there are
also linearly normal non-aCM  curves having minimal node postulation.
\end{rem}
\bigskip

\section{Curves on surfaces of small degree}\label{low degree}

In this section we study a class of curves having minimal node postulation,
with $n=3$,
$s\le 4$, $d> s(s-1)$ and $\alpha$ not minimal. This shows, in particular,
that the Halphen type bound for the genus in terms of $d$, $s$ and $e$
given in Corollary
\ref{bound for g}) is sharp in some non-trivial cases. We leave the easy
translation to
the reader.

\subsection{The case $s =2$}\label{s=2}

If $s = 2$ the inequality $d > s(s-1)$ is always satisfied, the minimal node
function $\Delta$
coincides with the lower bound given by Theorem \ref{pr} and is defined as
follows:
\begin{equation}\label{Delta 2}
\Delta(j) = \begin{cases}
0&\text{ if } \; j<0 \\
j+1 &\text{ if } \;  0\leq j \le \alpha-1\\
d-2-j&\text{ if } \;   \alpha \le j \le d-3\\
0&\text{ if } \;  j\ge d-2
\end{cases}.
\end {equation}
It follows that
\begin{equation} \label{delta Delta 2}
\delta(\Delta) = \frac 1 2\alpha(\alpha+1) + \frac 1 2
{(d-2-\alpha)}{(d-1-\alpha)}.
\end{equation}

\begin{prop} \label{curves on quadric 1} Assume $C$ lies on a quadric
surface $Q$. Then:
\begin{itemize}
\item[(a)] if $Q$ is smooth and $C$ is of type $(a,b)$ with $a \ge b$, then
$\alpha = a-1$;
\item[(b)] if $Q$ is a cone then
$$
\alpha = \begin{cases}
\frac{d}{2}-1&\text{ if } \; d \; \text{is even} \\
\frac{d-1}{2}&\text{ if } \; d \; \text{is odd}
\end{cases};
$$
\item[(c)] $C$ has minimal node postulation (\ref{delta2}).
\end{itemize}
\end{prop}

\begin{proof} (a) We have $e = b-2$ (e.g. by the adjunction formula) and the
conclusion follows by Remark \ref{adjoints}.

(b) If $d$ is even $C$ is a complete intersection of type $(2,\frac d 2)$ and
then $e = \frac d 2 - 2$, whence $\alpha = \frac d 2 -1$.

If $d$ is odd then $C$ is linked to a line $\ell$ in a complete
intersection of type
$(2,\frac{d+1} 2)$. Then the complete canonical series of $C$ is cut out by
the linear
system of all surfaces of degree $\frac{d+1} 2 - 2$ containing $\ell$,
outside the
divisor $\ell\cdot C$. It follows that $e = \frac{d+1} 2 - 3$, whence $\alpha =
\frac{d-1} 2$.

(c) Set $p := d-2-\alpha$. Then $p \ge 0$ and a straightforward calculation,
using (a), (b) and (\ref{delta Delta 2}) shows that
$$
\delta = \frac 1 2\alpha(\alpha+1) + \frac 1 2 p(p+1) = \delta(\Delta)
$$
The conclusion follows by Theorem \ref {pr} and an easy computation.
\end{proof}

The following Proposition gives a characterization of curves lying on a
quadric surfaces,
which completes Proposition \ref{curves on quadric 1}.

\begin{prop} \label {curves on quadric 2} Let $C \subseteq \pp^r$, $r \ge
3$, be a non degenerate curve. Then the following are equivalent:
\begin{itemize}
\item[(a)] $r = 3$, $C$ is special and is contained in a quadric;
\item[(b)] $\alpha \le d-4$ and $C$ has minimal node postulation
(\ref{Delta 2});
\item[(c)] $d\ge 6$ and $\partial h_N(d-4) = 2$;
\item[(d)] $C$ is special and $\partial h_N(d-4) = 2$.
\end{itemize}
\end{prop}

\begin{proof} (a) $\Rightarrow$ (b) follows from Remark \ref{adjoints}(v)
and Proposition \ref{curves on quadric 1}.

(b) $\Rightarrow$ (c). Since $\alpha \le d-4$ it is immediate that
$\partial h_N(d-4) = 2$.
Moreover $d-4 \ge \alpha \ge \partial h_N(d-4) = 2$, whence $d \ge 6$.

(c) $\Rightarrow$ (d). Assume that $C$ is non-special.
Then by Remark \ref{adjoints}(v) we have either $\alpha = d-3$ or $\alpha =
d-2$. In the
first case we have $d-3 = \partial h_N(d-4) = 2$, whence $d = 5$, a
contradiction. In the
second case we have, by Proposition \ref{basic}, $d-2 = \partial h_N(d-3) =
\partial
h_N(d-4) + 1 = 3$, whence $d = 5$, again a contradiction.

(d) $\Rightarrow$ (a). Since $C$ is special we have $\alpha \le d-4$.
Moreover by Theorem \ref {pr} we have $0 < \partial h_N(d-3) < \partial
h_N(d-4) = 2$,
whence $\partial h_N(d-3) = 1$. This implies $r = 3$.

Assume now $d \ge 7$. If $C$ is not contained in a quadric we can apply
Theorem \ref {p3bis} with $t = s = 3$ and $j = d-3$. This implies that
$\partial h_N(d-4)
\ge \partial h_N(d-3) + t-1 = 3$, a contradiction.

Finally assume $d = 6$. Then $d-3 = 3$ and we must have
$\partial h_N: \dots,  0,1,2,2,1,0,\dots$. Then by \cite {CGN}, Proposition
5.5, $C$ is
a complete intersection of type (2,3), and the conclusion follows also in
this case.
\end{proof}

As an application of the above Proposition we give the following
characterization of
the curves in $\pp^3$ having minimal $\alpha$.

\begin{prop}\label{alpha minimal} Assume $C\subseteq \pp^3$.
Then $\alpha = \lceil \frac{d-2}2 \rceil$ (that is $\alpha$ is minimal, see
Corollary
\ref{bound for alpha}) if and only if one (and only one) of the following
holds:
\begin{itemize}
\item[(i)] $(d,g) \in \{(3,0), (4,1), (5,1), (5,2)\}$ (in particular $C$ is
non-special);
\item[(ii)] $C$ is a special Castelnuovo curve;
\item[(iii)] $d = 9$ and $C$ is a complete intersection of two cubics;
\item[(iv)] $d = 7$ and $C$ is linked to a conic in a complete intersection
of type $(3,3)$.
\end{itemize}
\end{prop}

\begin{proof} Assume (i). Then $\alpha = d-2 = 1$  if $(d,g) = (3,0)$ and
$\alpha = d-3$ in the remaining cases. Then  $\alpha = \lceil \frac{d-2}2
\rceil$.

If (ii) holds then $\alpha$ is minimal by Proposition \ref{curves on
quadric 1} (see also remark \ref {quadric1'})

If either (iii) or (iv) hold then $\alpha$ is minimal by a direct calculation
(see e.g. Example \ref{CI} and Theorem \ref{Halphen}).

Assume now that $\alpha = \lceil \frac{d-2}2 \rceil$.

If $C$ is non-special then $d-3 \le \alpha \le d-2$.
If $\alpha = d-2$ then $d = 3$ and $g=0$.
If $\alpha = d-3$ we have either $d = 4$ or $d = 5$. In the first case
$\alpha = 1$ and
$\partial h_N: \dots 0,1,1,0\dots$, whence $g = 1$.
In the second case $\alpha = 2$ and either
$\partial h_N: \dots 0,1,2,1,0\dots$ (whence $g = 1$);
or $\partial h_N: \dots 0,1,2,2,0\dots$ (whence $g = 2$).

If $C$ is special and lies on a quadric then $C$ is a Castelnuovo curve by
Proposition \ref{curves on quadric 1}.

Finally assume that $C$ is special and does not lie on a quadric.
Then $d \ge 6$ and by \cite {CGN}, Proposition 5.5 $d$ must be odd. Then by
Proposition
\ref{HGP} (with $s=3$) we must have either $d = 9$ or $d = 7$. Moreover by
Proposition
\ref{curves on quadric 2} we must have
\begin{equation}\label{>2}
\partial h_N(d-4) > 2.
\end{equation}

Assume $d=9$, whence $\alpha = 4 < d-4$.
Then by Theorem \ref{pr} and (\ref{>2}) we must have $\partial h_N(d-4) =
3$ and
$\partial h_N(d-3) = 1$. It follows
$$
\partial h_N: \dots 0,1,2,3,4,4,3,1,0\dots,
$$
and then $C$ is a complete intersection of type $(3,3)$ by
Theorem \ref{Halphen} (with $d=9$, $s = 3$, $u = k =0$).

Assume finally that $d = 7$. Then an argument as the previous one shows that
$$
\partial h_N: \dots 0,1,2,3,3,1,0\dots.
$$
Then $C$ is linked to a conic in a complete intersection of type $(3,3)$, again
by Theorem \ref{Halphen} (with $d=7$, $s = 3$, $u = 0$, $k =2$).
\end{proof}

\begin{rem} \label{quadric1'}\rm
(i) Proposition \ref{curves on quadric 1}(c) shows that every curve in $\pp^3$
lying on a quadric, has minimal node postulation. This is true in
particular also by
non-Castelnuovo curves (that is curves of type $(a,b)$ on a smooth quadric,
with $|a-b| \ge 2$).

Moreover it is easy to check that every function (\ref{Delta 2}) occurs on
a smooth
quadric, provided $\alpha$ satisfies the obvious requirement $\alpha \ge
d-2-\alpha$
(just take curves of type $(\alpha +1, d-\alpha -1)$ on the quadric).
\smallskip

(ii) On the other hand Proposition \ref{curves on quadric 2} shows that
every {\it special} curve with minimal node postulation (\ref{Delta 2})
must lie on a
quadric, whence by Proposition \ref{alpha minimal} every special curve in
$\pp^3$
achieving both the bound of Theorem \ref{pr} and the lower bound for
$\alpha$ is
necessarily a Castelnuovo curve.
\smallskip

(iii) It follows easily by Remark \ref{adjoints}(v) that a non-special curve
$C \subseteq \pp^3$
has $\alpha\geq d-3$ hence it trivially attains the lower bound of Theorem
\ref{pr},
by Riemann-Roch and by \ref{adjoints} (iii). On the other hand such a $C$
does not
necessarily lie on a quadric surface (e.g. general rational curves of
degree $d \ge 5$
and elliptic curves of degree $d \ge 5$). Compare with Proposition
\ref{curves on quadric
2}.

(iv) We feel it would be interesting to generalize the above results to
higher-dimensional spaces, replacing quadric surfaces with surfaces of
minimal degree.
In particular we don't know if (ii) holds for such curves.
\end{rem}
\medskip

\subsection{The case $s = 3$}

The case where $\alpha$ is minimal and $C$ has the minimal node function is
yet
covered by Theorem \ref{Halphen} (compare condition (ii)). Hence we may
assume that
$\alpha$ is not minimal. It is easy to see that the minimal node function
can be written
as follows, for suitable integers $v\ge 0$ and $z\ge -1$:
\begin{equation} \label{Delta 3}
\Delta : \dots,\; 0,1,\;\dots,\; v+2z+1,\; 2z+1, \;\dots, \;3, 1, 0, \dots
\end{equation}
Here $d = v+3z+4$, $\alpha = 2z+v+1$ and if $v= 0$ the value $2z+1$ is
repeated twice,
that is $\Delta(\alpha) = \alpha$. Note also that in this case we have
$d-3-\alpha = z$,
hence any curve with minimal node postulation (\ref{delta2}) must have $e = z$.

From (\ref{Delta 3}) one computes immediately:
\begin{equation} \label{delta Delta 3}
\delta(\Delta) = \frac {(v+2z+1)(v+2z+2)} 2 + (z+1)^2.
\end{equation}
\medskip

\begin{prop} \label{curves on cubic} {\rm (a)} Let $S \subseteq \pp^3$ be a
smooth
cubic surface. Then for every pair of integers $(v,z)$ with $v \ge 0$ and
$z \ge -1$
there is a smooth connected curve $C \subseteq S$ having minimal node
postulation {\rm
(\ref{Delta 3})}.

{\rm (b)} Conversely any smooth curve $C \subseteq \pp^3$ having minimal node
postulation {\rm (\ref{Delta 3})} with $z \ge 2$ and $d> 12$ lies on a
cubic surface.
\end{prop}

\begin{proof} We use \cite{H}, Ch. V, \S 4 as a basic reference for the
geometry
of the cubic surface. We begin by showing that $S$ contains a smooth
rational curve $C'$
of degree $v+1$. Indeed (notation as in \cite {H}) it is sufficient to
consider a general
element of the linear system
$$
\Sigma :=
\begin{cases}
\left |\frac v 2\, l - (\frac v 2 - 1)e_1 \right |& {\rm if} \; v \; {\rm
is \; even}\\
&\\
\left |\frac {v+1} 2 \, l - (\frac {v-1} 2)e_1 - e_2\right |& {\rm if} \; v
\; {\rm is \; odd}
\end{cases}
$$

This shows our claim for $z = -1$. Assume now $z \ge 0$ and let $H$ be a
hyperplane
divisor on $S$. Then the linear system $|C'+(z+1)H|$ is very ample by
\cite{H}, V.4.12,
whence a general curve $C \in |C'+(z+1)H|$ is smooth and irreducible.

We claim that $C$ has the required properties. Indeed one computes easily
that $C$
has degree $v+3z+4 = d$ and genus $g = \frac{z(z+1)}2 +(z+1)(v+3)$, whence
$\delta =
\delta(\Delta)$ by (\ref{delta Delta 3}). Then by Theorem \ref{arit} $C$
has minimal
node postulation provided it has the right speciality, namely $e(C) = z$.

Now since $-H$ is a canonical divisor for $S$ we have that $(C'+zH)\cdot C$
is a
canonical divisor of $C$, whence $e(C) \ge z$. Then we need to show that
$(C'- H)\cdot
C$ is non-effective.

From the exact sequence
$$
0 \to \mathcal O_S(-H) \to \mathcal O_S(C'-H) \to \mathcal
O_{C'}((C'-H)\cdot C') \to 0,
$$
recalling that $C'$ is rational and that $(C'-H)\cdot C'$ is a canonical
divisor of $C'$,
we get $H^0(\mathcal O_S(C'-H)) = 0$. Hence from the exact sequence
$$
0 \to \mathcal O_S(-C+C'-H) \to \mathcal O_S(C'-H) \to \mathcal
O_C((C'-H)\cdot C) \to 0,
$$
we get $H^0(\mathcal O_C((C'-H)\cdot C)) = 0$, since
$H^1(\mathcal O_S(-C+C'-H)) =  H^1(\mathcal O_S(-(z+2)H) = 0$. This
concludes the proof
of (a).

To prove (b) let $C$ be a curve such that $\partial h_N = \Delta$ with $z
\ge 2$.
Then we have $\partial h_N(d-5) = 5$ and $\partial h_N(d-4) = 3$. On the
other hand,
since $d > 12$, we have:
$$
\binom {4-1}2 < d - (d-4) - 1 + \partial h_N(d-4) < d - \binom 4 2
$$
whence $h^0(\mathcal I_C(3)) > 0$ by Theorem \ref{p3bis}.
\end{proof}

\begin{rem}\label {no acm on cubic}\rm The curves with minimal node
postulation
we produce in the proof of Proposition \ref{curves on cubic} are, by
construction,
bilinked to the rational curve $C'$, whence they have, up to shift, the
same first
cohomology of $C'$ (see e.g. \cite {Mi}). In particular if $v+1= \deg(C')
\ge 4$ all our
curves are not aCM.

It would be interesting to know if there are aCM curves of the same degrees
and specialities, having minimal node postulation.
\end{rem}

\subsection {The case $s = 4$} Now we consider the case $s = 4$.
The case where $\alpha$ is minimal and $C$ has the minimal node function is
yet covered
by Theorem \ref{Halphen}(ii). Hence we may assume that $\alpha$ is not
minimal, namely
$\alpha > \lceil \frac {3d}4 \rceil - 3$. We shall also assume that $\alpha
\le d-5$,
otherwise the meaningful part of the node function, from $\alpha$ to $d-3$,
is too short
and everything trivializes. Then the minimal node function can be written
as follows, for
suitable integers $v \ge 1$ and $z \ge 1$:
\begin{equation} \label{Delta 4}
\Delta : \dots,\; 0,1,\;\dots,\; v+3z+3,\; 3z+3, \;\dots, \; 6, \;3, 1, 0,
\dots
\end{equation}
Here $d = v+4z+7$ and $\alpha = 3z+v+3$. Note that $d-3-\alpha = z+1$,
whence if a curve has minimal node postulation then necessarily $e = z+1$.

From (\ref{Delta 4}) one computes immediately:
\begin{equation} \label{delta Delta 4}
\delta(\Delta) = \frac {(v+3z+3)(v+3z+4)} 2 + \frac {3(z+1)(z+2)} 2 +1.
\end{equation}
\medskip

\begin{prop} \label{curves on quartic} {\rm (a)} For every pair of integers
$(v,z)$
with $v \ge 1$ and $z \ge 1$ there are a smooth quartic surface $S
\subseteq \pp^3$
and a smooth connected curve $C \subseteq S$ having minimal node postulation
{\rm (\ref{Delta 4})}.

{\rm (b)} Conversely any smooth curve $C \subseteq \pp^3$ having minimal node
postulation {\rm (\ref{Delta 4})} with $z \ge 2$ and $d > 20$ lies on a
quartic surface.
\end{prop}

\begin{proof} The proof is similar to the previous one. By \cite{M} there
are a
smooth quartic surface $S$ and a smooth connected curve $C' \subseteq S$
with the
following properties:
\begin{itemize}
\item[(i)] $C'$ has genus $v$ and degree $v+3$;
\item[(ii)] $\Pic(S)$ is freely generated by the classes of $C'$ and $H$,
where $H$ is a hyperplane divisor.
\end{itemize}

Since $K_S \sim 0$ we have $D^2 = 2p_a(D) - 2$ for every divisor $D$.
In particular ${C'}^2 = 2v-2 \ge 0$, so that $|C'|$ has no base points,
whence $|C' +
(z+1)H|$ is very ample. Then a general $C \in |C' + (z+1)H|$ is smooth and
connected and
we will show that $C$ has the required properties.

First of all one computes $\deg(C) = v + 4z + 7$ and $g = v + (z+1)(v+3) +
2(z+1)^2$,
whence $\delta = \delta(\Delta)$ by (\ref{delta Delta 4}).

It remains to show that $C$ has the required speciality, namely $e(C) =
z+1$. ùIt is
clear that $e(C) \ge z+1$ and then we have to show that the linear system
$|C' - H|$ is
non-effective (see proof of Proposition \ref{curves on cubic}).

We assume then $|C' - H|$ effective and we seek for a contradiction.

Since $(C' - H)^2 = - 4$ this linear system must have a fixed integral
curve $E$
with negative self-intersection. It follows that $E$ is rational and $E^2 =
- 2$, whence,
in particular,
$|E| \neq |C' - H|$. Then  if $p := \deg (E)$ we have $p \le v-2$.
On the other hand, by (ii)
above, there are uniquely determined integers $x$ and $y$ such that $E \sim
xC' + yH$,
whence:
\begin{equation}\label {deg E} (v+3)x + 4y = p
\end{equation}
and
\begin{equation}\label {E2}
2(v-1)x^2 + 2(v+3)xy + 4y^2 = -2
\end{equation}

Squaring (\ref{deg E}) and subtracting (\ref{E2}) multiplied by $4$ we get
$$
x^2 = \frac {p^2+8}{(v+3)^2-8(v-1)},
$$
whence $x^2 < 1$ by a simple computation, remembering that $p \le v-2$.
This is a
contradiction and the proof of (a) is complete.

Now we prove (b).  Let $\Delta$ be as in (\ref{Delta 4}) with $z \ge 2$ and
let $C$ be a curve such that $\partial h_N = \Delta$. Then  we have
$\partial h_N(d-6) =
9$ and $\partial h_N(d-5) = 6$. On the other hand, since $d > 20$, we have:
$$
\binom {5-1}2 < d - (d-5) - 1 + \partial h_N(d-5) < d - \binom 5 2
$$
whence $h^0(\mathcal I_C(4)) > 0$ by Theorem \ref{p3bis}.
\end{proof}

The next example shows that there are smooth quartic surfaces, other than
the Mori
surfaces, containing curves with minimal nodal postulation.

\begin{ex}\label{noaCM} \rm Let  $C$ be directly linked to a pair $Y$ of
skew lines
by surfaces of degree $4,5$. One computes $d=18$, $s=4$ (via the mapping
cone) and one
gets $\alpha$ reminding that $\omega_{C}$ is cut by surfaces of degree
$5+4-4=5$ passing
through $Y$ (outside $Y$) whence $e = 3$, because $Y$ lies on a quadric. It
follows that
$\alpha = 12$. Moreover
$h^1(\mathcal I_C(1)) = h^1(\mathcal I_Y(4)) = 0$ by liaison, whence
$h^0(\mathcal O_C(1)) = 4$ and $n = 3$ \par\noindent

The minimal node function $\Delta$ is:
$$
 ...,0,1,2,3,4,5,6,7,8,9,10,11,12,9,6,3,1,0,...
$$
and its sum is $\delta(\Delta) = 97$. On the other hand by liaison we have
$g = 39$
whence $\delta = 97$. Then $C$ has minimal nodal postulation by Corollary
\ref{bound for
g}.

Observe also that in this case $\alpha$ is not minimal (the least value
being $11$),
in agreement with Theorem \ref{Halphen}.

Note also that the quartic surface cannot be one of the surfaces considered
in the proof
of Proposition \ref{curves on quartic}, because the classes of the two skew
lines and of a
hyperplane divisor are easily seen to be linearly independent.
\end {ex}

\begin{rem} \rm If $v > 3$ the curves we produce in the proof of Proposition
\ref{curves on quartic} are not aCM. Indeed the starting curve $C'$ of
degree $v+3$ and
genus $v$ has non trivial $H^1(\mathcal I_{C'}(2))$ for $v > 3$, and then
we can use a
bilinkage argument as in Remark \ref{no acm on cubic}.
\end{rem}

\end{document}